\documentclass{tpms-l}

\newtheorem{theorem}{Theorem}[section]

\theoremstyle{definition}

\newtheorem{example}[theorem]{Example}

\theoremstyle{remark}
\newtheorem{remark}[theorem]{Remark}

\numberwithin{equation}{section}

\usepackage[utf8]{inputenc}
\usepackage[T1]{fontenc}
\usepackage[british]{babel}
\usepackage{amsmath}

\usepackage{graphicx}
\usepackage{xcolor}
\usepackage{xurl}
\usepackage{caption}

\begin{document}

\title{On Spectral Theory of Random Fields in the Ball}

\author{Nikolai Leonenko}
\address{School of Mathematics, Cardiff University, Senghennydd Road, Cardiff CF24 4AG, UK}
\email{LeonenkoN@Cardiff.ac.uk}

\author{Anatoliy Malyarenko}
\address{Division of Mathematics and Physics, Mälardalen University, 721 23 Västerås, Sweden}
\email{anatoliy.malyarenko@mdh.se}

\author{Andriy~Olenko}
\address{Department of Mathematics and Statistics, La Trobe University,  Melbourne, VIC 3086, Australia}
\email{A.Olenko@latrobe.edu.au}


\date{\today}

\dedicatory{The paper is dedicated to the 90th birthday of Professor Myhailo Yosypovych Yadrenko (1932--2004).}


\begin{abstract}
The paper investigates random fields in the ball. It studies three types of such fields: restrictions of scalar random fields in the ball to the sphere, spin, and vector random fields. The review of the existing results and new spectral theory for each of these classes of random fields are given.  Examples of applications to classical and new models of these three types are presented. In particular,  the Mat\'{e}rn model is used for illustrative examples. The derived  spectral representations can be utilised to further study theoretical properties of such fields and to simulate their realisations.  The obtained results can also find various applications for modelling and investigating ball data in cosmology, geosciences and embryology.
\end{abstract}

\maketitle

\section{Introduction}

Recent years have witnessed an enormous amount of attention to investigating spherical random fields. The theoretical interest (see, for example, \cite{ Malyarenko:2012, Marinucci:2011, Yadrenko_book} and the references therein) is strongly influenced by studies of random fields on manifolds, as the sphere is one of the simplest manifolds.  The empirical motivation comes from cosmology, earth science and embryology, just to name a few (see, for instance, \cite{Michel, Pike, Porcu, Weinberg:2008}). The main approaches and tools in such investigations are based on the spectral theory of spherical random fields. Professor Yadrenko was one of pioneering researchers and leading figures in developing this theory. Later on, it was demonstrated that the behaviour of the power angular spectrum determines various properties of these fields and evolutions of their spatio-temporal counterparts, see \cite{Broadbridge:2019,  Broadbridge:2020, Lang:2015}. However, the known results about spherical fields are not directly translatable to the random fields defined in the ball. Therefore, most of the spectral theory for different classes of such fields should be developed independently.

One of main applied motivations for developing the spectral theory of random fields in the ball comes from cosmological research. The future European Space Agency mission Euclid and Cosmic Microwave Background Stage 4 (CMB-S4) project supported by the US Department of Energy Office of Science and the National Science Foundation are planned to collect and analyse cosmological data in a ball of radius about $10$~billion light years. From the mathematical point of view, these missions will sample values of several scalar, spin and tensor random fields defined in the ball. It requires further development of stochastic models and statistical tools for such fields.

Deterministic spin fields on the sphere were introduced by \cite{MR0047664}. They became well known to physicists after the seminal paper \cite{MR194172}. Random spin fields on the sphere appeared in \cite{PhysRevD.55.1830} as a technical tool for analysing a full-sky polarisation map of the cosmic microwave background. This problem was also  independently studied in \cite{PhysRevD.55.7368} by using tensor random fields on the sphere.  A comprehensive review of deterministic spin and tensor fields on the sphere can be found in \cite{MR569166}.

In stochastic settings, the rigorous mathematical theory of spin random fields on the sphere was proposed by \cite{MR3170229}, \cite{MR2737761}, and \cite{MR2884225}. This theory works well for studies of the current cosmic microwave background radiation data collected on the sphere. However, modelling and statistical analysis of data from the Euclid and CMB-S4 surveys requires a generalisation of the above theory to random fields in the ball. First steps of such generalisation were proposed by \cite{Anatoliy_Malyarenko_2015}. One of main ideas, that was originally suggested by M.~Yadrenko in \cite{MR0164376}, is outlined in Section~\ref{sec:isotr}.

This paper studies three main classes of random fields in the ball: restrictions of scalar random fields in the ball to the sphere, isotropic spin,  and vector  random fields. It presents some existing in the literature results and develops new representations for those cases that were not covered before. It suggests a unified approach and notations in the spectral theory of random fields in the ball. The results could be useful for further studying and comparing the three classes mentioned above. Several examples of applications to classical and new models provide explicit spectral representations, which can be used in spatial statistics. To the best of our knowledge the explicit expressions for spectral coefficients of the Mat\'{e}rn model are also new.  All coefficients in the derived theoretical representations are easily computable and can be utilised in numerical applications.

The structure of the paper is as follows. Section~\ref{sec:isotr} presents main definitions and results about isotropic random fields  that are obtained via restrictions of random fields in the ball to the sphere. Results about spin  random field on the sphere are given in Section~\ref{sec:sphere}. The spectral theory of spin  random field in the ball is developed in Section~\ref{sec:ball}. Section~\ref{sec:multi} studies spectral properties of vector $\rho$-stationary random fields.  Finally, the conclusions and some future research directions are presented in~Section~\ref{sec:conclusion}.

All numerical examples were produced by using the software Maple version 2021.0. This software was also used to verify some theoretical computations.
A reproducible version of the code in this paper is available in the folder "Research materials" from the website~\url{https://sites.google.com/site/olenkoandriy/}.

\section{Spherical isotropic random fields as restrictions of fields in the ball}\label{sec:isotr}

Let us denote the centered ball of radius $r_0>0$ by
\[
\mathbb{B}(r_0)=\{\,\mathbf{x}\in\mathbb{R}^3\colon\|\mathbf{x}\|\leq r_0\,\},
\]
where $\|\cdot\|$ denotes the Euclidean norm in ${\mathbb{R}}^3$.

Let $T(\mathbf{x})$, $\mathbf{x}\in\mathbb{R}^3\ (\mbox{or}\ \mathbf{x}\in  \mathbb{B}(r_0)),$ be a random field. In other words, there is a probability space $(\Omega,\mathfrak{F},\mathsf{P})$ and a function $T\colon\mathbb{R}^3\times\Omega\to\mathbb{C}$ such that for any fixed $\mathbf{x}\in\mathbb{R}^3$ the function $T(\mathbf{x},\omega)$ is a complex-valued random variable. Assume that the random field $T(\mathbf{x})$ is second-order, that is, $\mathsf{E}[|T(\mathbf{x})|^2]<\infty,$  and mean-square continuous, that is,
\[
\lim_{\mathbf{y}\in\mathbb{R}^3:\|\mathbf{y}-\mathbf{x}\|\to 0}\mathsf{E}[|T(\mathbf{y})-T(\mathbf{x})|^2]=0,
\quad \mbox{for all} \quad \mathbf{x}\in\mathbb{R}^3.
\]
 
Let $\langle T(\mathbf{x})\rangle=\mathsf{E}[T(\mathbf{x})]$ be the one-point correlation function of the random field $T(\mathbf{x})$, and let
\[
\langle T(\mathbf{x}),T(\mathbf{y})\rangle=\mathsf{E}[\overline{(T(\mathbf{x})-\langle T(\mathbf{x})\rangle)}(T(\mathbf{y})-\langle T(\mathbf{y})\rangle)]
\]
be its two-point correlation function. Let $G=\mathrm{SO}(3)$ be the rotation  group in $\mathbb{R}^3,$ that is, the group of orthogonal $3\times 3$ matrices with a unit determinant.

Call the field $T(\mathbf{x})$ {isotropic} if its one-point correlation function is constant, while its two-point correlation function is rotation-invariant:
\[
\langle T(g\mathbf{x}),T(g\mathbf{y})\rangle=\langle T(\mathbf{x}),T(\mathbf{y})\rangle,\qquad g\in\mathrm{SO}(3).
\]

Without loss of generality,  this paper assumes that $\langle T(\mathbf{x})\rangle=0.$

How to describe the class of all possible two-point correlation functions of isotropic random fields? Consider the restriction of the field $T(\mathbf{x})$ to  $\mathbb{S}^{2}(r),$ which denotes the centred sphere of radius $r>0$ in $\mathbb{R}^3$. To avoid introducing new notations $T(r,\theta,\varphi)$ will be used for $(r,\theta,\varphi)$ which are the spherical coordinates of $\mathbf{x}.$ The two-point correlation function of the above restriction is rotation-invariant and depends only on the angle between two points. Thus, the restriction is an isotropic field on the sphere. Such fields were completely described by \cite{obukhov} and have the form
\begin{equation}\label{eq:1}
T(r,\theta,\varphi)=\sum_{\ell=0}^{\infty}\sum_{m=-\ell}^{\ell}a_{\ell m}(r)Y_{\ell m}(\theta,\varphi),
\end{equation}
where $(r,\theta,\varphi),$ $r>0,$ $\theta \in[0, \pi],$ $\varphi \in [0, 2\pi)$ are the spherical coordinates of a point $\mathbf{x}\in\mathbb{R}^3$, $\{Y_{\ell m}(\theta,\varphi), \ell  \in\mathbb N_0, m=-\ell,...,\ell \} $ with $\mathbb N_0= \mathbb N \cup \{0\},$ are the spherical harmonics, and  $a_{\ell m}(r)$ are finite variance random variables
\[
a_{\ell m}(r)=\int_{0}^{\pi}\int_{0}^{2\pi}T(r,\theta,\varphi)
\overline{Y_{\ell m}(\theta,\varphi)}\sin\theta\,\mathrm{d}\theta
\,\mathrm{d}\varphi
\]
that satisfy the conditions
\begin{equation}\label{eq:2}
\begin{aligned}
\mathsf{E}[a_{\ell m}]&=0,\\
\mathsf{E}[a_{\ell m}(r)\overline{a_{\ell'm'}(r)}]&=\delta_{\ell\ell'}\delta_{mm'}C_{\ell}(r)
\end{aligned}
\end{equation}
for all $\ell,\ell'\in \mathbb N_0,$ $ m=-\ell,...,\ell,$  $m'=-\ell',...,\ell'.$

For each $r>0,$ the sequence  $\{C_{\ell}(r),\ell\in \mathbb N_0\}$ of non-negative numbers satisfies the condition
\[
\sum_{\ell=0}^{\infty}(2\ell+1)C_{\ell}(r)<\infty.
\]

As a function of the variable $r$, $a_{\ell m}(r)$ is a stochastic process. Denote
\[
C_{\ell}(r_1,r_2)=\mathsf{E}[a_{\ell m}(r_1)\overline{a_{\ell'm'}(r_2)}].
\]
It follows from \eqref{eq:1} and \eqref{eq:2} that
\[
\mathsf{E}[T(\mathbf{x}_1)\overline{T(\mathbf{x}_2)}]=\sum_{\ell=0}^{\infty}
\sum_{m=-\ell}^{\ell}Y_{\ell m}(\theta_1,\varphi_1)\overline{Y_{\ell m}(\theta_2,\varphi_2)}C_{\ell}(r_1,r_2).
\]
The addition theorem for spherical harmonics implies that
\[
\mathsf{E}[T(\mathbf{x}_1)\overline{T(\mathbf{x}_2)}]=\frac{1}{4\pi}
\sum_{\ell=0}^{\infty}(2\ell+1)C_{\ell}(r_1,r_2)P_{\ell}(\cos\gamma),
\]
where $\gamma$ is the angle between the vectors $\mathbf{x}_1$ and $\mathbf{x}_2$ and $\{P_{\ell}(\cdot),\ell\in \mathbb N_0\}$ are the Legendre polynomials.

If $T(\mathbf{x}),\mathbf{x}\in \mathbb{R}^{3},$ is a homogeneous and
isotropic random field, then its covariance function has the following spectral representation
\begin{equation*}
\left\langle T(\mathbf{x}),T(\mathbf{y})\right\rangle =\int_{0}^{\infty}\frac{\sin
(\lambda \left\Vert \mathbf{y}-\mathbf{x}\right\Vert )}{\lambda \left\Vert \mathbf{y}-\mathbf{x}\right\Vert }%
d\mu (\lambda ),\qquad\mathbf{x},\mathbf{y}\in \mathbb{R}^{3},
\end{equation*}%
where $\mu $ is the finite measure.

Therefore, for the random field (\ref{eq:1}) on the sphere  $\mathbb{S}^{2}(r)$ it holds
\begin{equation*}
C_{\ell}(r)=2\pi\int_{0}^{\infty}\frac{J_{\ell+\frac{1}{2}%
}^{2}(\lambda r)}{\lambda r}d\mu (\lambda ),\quad \ell\in \mathbb N_0,
\end{equation*}%
where $J_{{}\nu }\left( z\right) $ is the Bessel function of the first kind
of order $\nu.$

In this case%
\begin{equation*}
\left\langle T(\mathbf{x}),T(\mathbf{y})\right\rangle =B(\left\Vert \mathbf{y}-\mathbf{x}\right\Vert )=\int_{0}^{\infty }\frac{\sin \left(2r\lambda
\sin \left(\frac{\gamma }{2}\right)\right)}{2r\lambda \sin  \left(\frac{\gamma }{2}\right)}d\mu (\lambda ),
\end{equation*}%
where the Euclidean distance  $\left\Vert
\mathbf{y}-\mathbf{x}\right\Vert,$ called also the
\emph{chordal distance}, between two points on a sphere $\mathbf{x},\mathbf{y}\in  \mathbb{S}^{2}(r)\subset \mathbb{R}^{3},$  can be expressed in terms of the great circle (also known
as geodesic or spherical) distance as follows:%
\begin{equation*}
\left\Vert \mathbf{y}-\mathbf{x}\right\Vert =2r\sin \left(\frac{\gamma }{2}\right),
\end{equation*}%
where $\gamma =\mathbb{\gamma (}\mathbf{x},\mathbf{y})=\arccos \left\langle \mathbf{x},\mathbf{y}\right\rangle $
 and $\left\langle \mathbf{x},\mathbf{y}\right\rangle $ denotes the usual inner product on $\mathbb{R}^{3}$.

\begin{example}[Mat\'{e}rn covariance function]\label{ex1}

Consider a covariance function of a scalar random field $T(\mathbf{x}),\mathbf{x}\in
\mathbb{R}^{3},$ of the form
\begin{equation}\label{mat1}
\left\langle T(\mathbf{x}),T(\mathbf{y})\right\rangle =\sigma ^{2}\,\frac{%
2^{1-\nu }}{\Gamma \left( \nu \right) }\left( a\left\Vert
\mathbf{y}-\mathbf{x}\right\Vert \right) ^{{}\nu }K_{{}\nu }\left( a\left\Vert \mathbf{y}-\mathbf{x}\right\Vert
\right) ,\quad
\end{equation}%
where $\sigma ^{2}>0,$ $a>0,$ $\nu >0,$ and $K_{{}\nu }\left( \cdot\right) $ is the modified Bessel function of the second kind of order $\nu .$ Here, the parameter $\nu $
controls the degree of differentiability of the random field,  $\sigma $
is field's variance and the parameter $a$ is a scale parameter.

 The corresponding
isotropic spectral density is
\begin{equation*}
f\left( \lambda\right) =\frac{\sigma
^{2}\Gamma \left( \nu +\frac{3}{2}\right) a^{2\nu }}{\pi ^{3/2}\Gamma(\nu)\left(
a^{2}+\lambda^{2}\right) ^{\nu +\frac{3}{2}}},\quad \lambda\geq 0.
\end{equation*}

The restriction of an homogeneous and isotropic Mat\'{e}rn random
field to the sphere $\mathbb{S}^{2}(r)$ is an isotropic field on this sphere with the covariance structure
\begin{equation*}
\begin{aligned}
\left\langle T(\mathbf{x}),T(\mathbf{y})\right\rangle = B\left(2\sin \left(\frac{\gamma }{2}\right)\right)=\frac{2^{1-\nu }\sigma ^{2}}{\Gamma
\left( \nu \right) }\left( a^2\sin \left(\frac{\gamma }{2}\right)\right) ^{{}\nu }K_{{}\nu
}\left( a^2\sin \left(\frac{\gamma }{2}\right)\right) ,
\end{aligned}
\end{equation*}%
while the application of the formula (12) from \cite[\S1.6]{Yadrenko_book} results in the angular spectrum of the form
\begin{equation*}
C_{\ell}(r)=4\pi^{3/2}\sigma ^{2}\frac{\Gamma \left( \nu +\frac{3}{2}\right) a^{2\nu }}{\Gamma(\nu)\,r}\int_{0}^{\infty }J_{\ell+\frac{1}{2}}^{2}(r\lambda)\lambda\left(
a^{2}+\lambda^{2}\right) ^{-(\nu +\frac{3}{2})}d\lambda,\ \ell\in \mathbb N_0.
\end{equation*}

To calculate this integral, once can use \cite[Equation~2.12.32.10]{MR950173} and obtain
\[
\begin{aligned}
	C_{\ell}(r)&=2\pi^{3/2}\sigma ^{2} a^{2\nu }\left(\frac{\nu\,\Gamma(\ell-\nu)}{\sqrt{\pi}\,\Gamma(\ell+\nu+2)}
	{}_1F_2(\nu+1;\nu-\ell+1,\nu+\ell+2;a^2r^2)\,r^{2\nu}\right.\\
	&\quad+\left.
	\frac{\Gamma(\nu-\ell)\,a^{2\ell-2\nu}}{2^{2\ell+1}\Gamma(\nu)\Gamma(\ell+3/2)}
	{}_1F_2(\ell+1;\ell-\nu+1,2\ell+2;a^2r^2)\,r^{2\ell}\right),
\end{aligned}
\]
where ${}_1F_2$ is the generalised hypergeometric function. For zero and negative integer values of $\ell-\nu$ or $\nu-\ell$ the above expression is  interpreted as its limit when $\nu$ approaches~$\ell.$ The limit is finite due to the asymptotic behaviour of the generalised hypergeometric function ${}_1F_2(\cdot).$

For specific values of the parameters the expressions above can be simplified to the forms that can be easily used in computations. For example, for $a=10,$ $\sigma^2=1,$ and $\nu=1/2$ one obtains

\[ B\left(2\sin \left(\frac{\gamma }{2}\right)\right)=\exp\left(-20 \sin \left(\frac{\gamma }{2}\right) \right), \quad f\left( \lambda\right)=\frac{10}{\pi ^{2}\left(
100+\lambda^{2}\right)^2},
\]
\[
\begin{aligned} C_{\ell}(r)&= \frac{{\pi}}{5 r}\left(10rI_{\ell +\frac{1}{2}}\! \left(10 r \right)K_{\ell +\frac{3}{2}}\! \left(10 r \right)-10rK_{\ell +\frac{1}{2}}\! \left(10 r \right)I_{\ell +\frac{3}{2}}\! \left(10 r \right)\right.\\
&\quad -\left.
(2\ell+1)I_{\ell +\frac{1}{2}}\! \left(10 r \right) K_{\ell +\frac{1}{2}}\! \left(10 r \right)\right),
\end{aligned}
\]
where $I_{{}l }\left( \cdot\right) $ is the modified Bessel function of the first kind of order $l.$

The plot of the covariance function (\ref{mat1}) is shown in Figure~\ref{fig1}. To produce this plot the values $\mathbf{x}=\mathbf{0}$ and $\mathbf{y}=(y_1,y_2,y_3)\in \mathbb B_0(r_0)$ with $y_3=0$ were chosen. The horizontal coordinates are $(y_1,y_2),$ while the vertical one represents the values of  $\left\langle T(\mathbf{0}),T(\mathbf{y})\right\rangle.$

Plots of first few coefficients $C_{\ell}(r)$ of the corresponding angular power spectrum on the interval $r\in[0,1]$ are given in Figure~\ref{fig2}.

\begin{figure}[!htb]
 \captionsetup{width=0.45\textwidth}
	\centering
	\begin{minipage}{0.48\linewidth}
			\includegraphics[width=1\linewidth,height=1\textwidth,trim =5mm 0 15mm 0 ,clip]{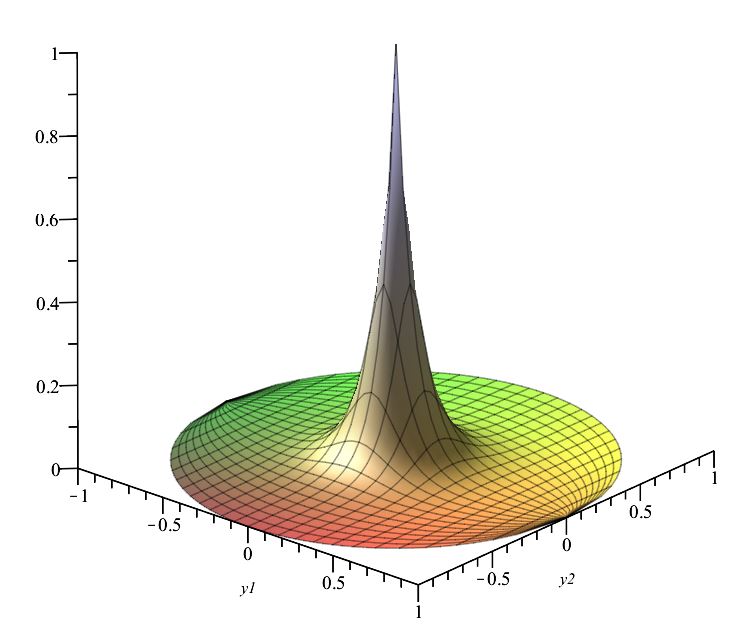}
			\caption{Mat\'{e}rn covariance function for $a=10,$ $\sigma^2=1,$ and  $\nu=1/2$.}\label{fig1}
	\end{minipage}
	\hfill
	\begin{minipage}{0.48\textwidth}
			\includegraphics[width=1\textwidth,height=1\textwidth,trim =5mm 0 5mm 0 ,clip]{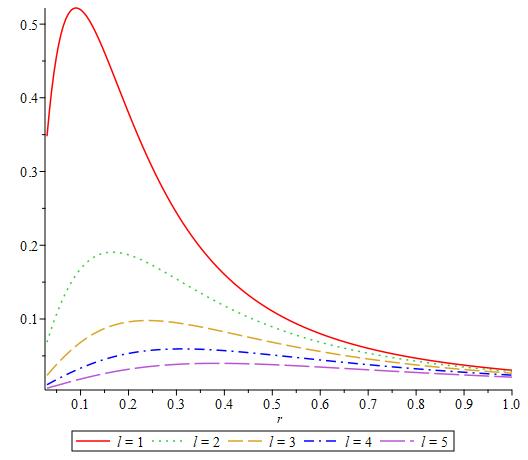}
			\captionof{figure}{ $C_{\ell}(r)$ for $a=10,$ $\sigma^2=1,$ and  $\nu=1/2$.}\label{fig2}
	\end{minipage}
	\end{figure}

\end{example}

\section{Spin random fields on the sphere}\label{sec:sphere}

To define spin and tensor random fields in the ball, the opposite direction is used. Let $T(\mathbf{x})$ be a random field defined in the centered ball
$\mathbb{B}(r_0).$
Call the field $T(\mathbf{x})$ spin or tensor, if for any $r\in(0,r_0]$ the restriction of the field to the centred sphere of radius~$r$ is a spin or  tensor random field on this sphere $\mathbb{S}^{2}(r).$ Starting from  results about spin or tensor random fields on the sphere, we will construct the spectral theory of spin or tensor random fields in $\mathbb{B}(r_0).$

There are two different approaches to deterministic spin fields on a manifold, see~\cite{MR2384313}. The first one requires introducing the so-called principal bundles of spin frames and will not be introduced here. The second one is as follows.

Let $(E,\pi,M)$ be a vector bundle over a manifold $M$. In particular, $\pi\colon E\to M$ and there is an open covering $\{U_{\alpha}\}$ of $M$, a finite-dimensional linear space $L$, and the one-to-one maps $\varphi_{\alpha}\colon\pi^{-1}(U_{\alpha})\to U_{\alpha}\times L$ such that for all $x\in U_{\alpha}$ the set $\pi^{-1}(x)$ is a copy of $L$, and the overlaps $\varphi_{\alpha}\circ\varphi^{-1}_{\beta}$ map a point $(x,v)\in(U_{\alpha}\cap U_{\beta})\times L$ to a point $(x,g_{\alpha\beta}v)$ for some suitable change-of-coordinates invertible linear operators $g_{\alpha\beta}(x)$.

Various conditions on $M$ can be formulated in terms of the functions $g_{\alpha\beta}(x)$. For example, $M$ is orientable if and only if there is such an open covering $\{U_{\alpha}\}$ of $M$, that the above functions take values in the connected component of unity of the group $\mathrm{GL}(L)$ of invertible linear operators on $L$. $M$ is orientable and Riemannian if and only if $L$ is a real linear space and the functions $g_{\alpha\beta}(x)$ take values in the group $\mathrm{SO}(L)$ of orthogonal linear operators with unit determinant for a suitable covering. Finally, $M$ is \emph{spin} if the space $L$ carries a special representation of the so-called \emph{spin group} that covers the group $\mathrm{SO}(L)$ twice. For details, see \cite{MR1031992}. Both the sphere and the ball are spin manifolds, and spin random fields can be properly defined on them.

We remind the results of the general theory of spin random fields on the sphere, see~\cite{MR3170229}, \cite{MR2737761} and \cite{MR2884225}. Let $s$ be an integer, and let $K=\mathrm{SO}(2)$ be the group of rotations of the three-dimensional space around the $z$-axis. Each element of  $\mathrm{SO}(2)$ can be given in the form
\[
k_{\varphi}=
\begin{pmatrix}
  \cos\varphi & \sin\varphi & 0 \\
  -\sin\varphi & \cos\varphi & 0 \\
  0 & 0 & 1
\end{pmatrix}
,\qquad 0\leq\varphi<2\pi.
\]
The correspondence that maps this element to the $1\times 1$ unitary matrix $\mathrm{e}^{-\mathrm{i}s\varphi}$ is an irreducible unitary representation of the group $\mathrm{SO}(2)$. Consider the Cartesian product $\mathrm{SO}(3)\times\mathbb{C}$. Call two elements $(g_1,z_1)$ and $(g_2,z_2)$ in $\mathrm{SO}(3)\times\mathbb{C}$ equivalent if there exists a $\varphi$ such that $(g_2,z_2)=(g_1k_{\varphi},\mathrm{e}^{\mathrm{i}s\varphi}z_1)$. Call the set of equivalence classes $E_s$.

Let $p\colon\mathrm{SO}(3)\times\mathbb{C}\to E_s$ be the correspondence that maps an element $(g,z)$ to its equivalence class. Equip $E_s$ with the quotient topology, that is, a set $A\subseteq E_s$ is open if and only if its inverse image $p^{-1}(A)$ is open in $\mathrm{SO}(3)\times\mathbb{C}$. Consider the mapping $\pi$ that maps an element $(g,z)$ to the left coset $gK$. All elements of the same equivalence class have the same image under $\pi$, so one can consider $E_s$ as the domain of $\pi$. The image of~$E_s$ under $\pi$ is the set $G/K$ of all left cosets, which is the centred unit sphere $\mathbb{S}^{2}\subset\mathbb{R}^3$. The triple $(E_s,\pi,\mathbb{S}^2)$ is a line bundle over $\mathbb{S}^2$.

A mapping $f\colon \mathbb{S}^2\to E_s$ is called a section of the line bundle $(E_s,\pi,\mathbb{S}^2)$ if $f(\mathbf{x})\in\pi^{-1}(\mathbf{x})$ for all $\mathbf{x}\in \mathbb{S}^2$. Let $\mu$ be the Lebesque measure on $\mathbb{S}^2$. Let $L^2(E_s)$ be the set of $\mu$-equivalence classes of all sections $f$ with
\[
\int_{\mathbb{S}^2}|f(\mathbf{x})|^2\,\mathrm{d}\mu(\mathbf{x})<\infty.
\]
Equation
\[
U(g)f(\mathbf{x})=f(g^{-1}\mathbf{x}),\qquad g\in\mathrm{SO(3)},
\]
defines a unitary representation of the group $\mathrm{SO}(3)$ in the complex Hilbert space $L^2(E_s)$.

The irreducible unitary representations of the group $\mathrm{SO}(3)$ are enumerated by non-negative integers $\ell$ (this is the traditional notation of the angular momentum in quantum mechanics). Let $(\alpha,\beta,\gamma)$ be the Euler angles of a rotation $g\in\mathrm{SO}(3)$. There are many different conventions in the literature, see \cite{MR2884225} for a survey. Here and in what follows we adopt conventions from \cite{durrer2008cosmic}. In particular, the first rotation is by angle $\gamma$ around the $z$-axis, then a rotation by angle $\beta$ around the $y$-axis and finally a rotation by angle $\alpha$ around the new $z$-axis.

Let $D^{(\ell)}_{m,n}(\alpha,\beta,\gamma)$ be the matrix entries of the $\ell$th irreducible unitary representation in the basis described in \cite[p.~344]{durrer2008cosmic}. The sections of the line bundle $(E_s,\pi,\mathbb{S}^2)$ defined in the local chart of spherical coordinates by
\[
{}_sY_{\ell m}(\theta,\varphi)=\sqrt{\frac{2\ell+1}{4\pi}}\,
\overline{D^{\ell}_{m,-s}(\varphi,\theta,0)},\quad  \ell  \ge s,\  m=-\ell,...,\ell,
\]
are called the \emph{spin weighted spherical harmonics}. They are defined for $\ell\ge s$ and $|m|\leq\ell$ and form an orthonormal basis in the space $L^2(E_s)$:
\[
\int_{0}^{\pi}\int_{0}^{2\pi}{}_sY_{\ell m}(\theta,\varphi)
\overline{{}_sY_{\ell'm'}(\theta,\varphi)}\sin\theta\,\mathrm{d}\theta
\,\mathrm{d}\varphi=\delta_{\ell\ell'}\delta_{mm'}.
\]
Note that ${}_sY_{\ell m}(\theta,\varphi)$ are functions on $S^2$ if and only if $s=0$. Otherwise, they are sections of a nontrivial bundle $(E_s,\pi,\mathbb{S}^2)$ that cannot be represented as a Cartesian product $S^2\times\mathbb{C}$.

A random section ${}_sT(\mathbf{x})$ of the line bundle $(E_s,\pi,\mathbb{S}^2)$ is called an isotropic spin $s$ random field if for all $\mathbf{x}$, $\mathbf{y}\in \mathbb{S}^2$, and for all $g\in\mathrm{SO}(3)$ it holds
\[
\begin{aligned}
\langle {}_sT(g\mathbf{x})\rangle&=\langle {}_sT(\mathbf{x})\rangle,\\
\langle {}_sT(g\mathbf{x}),{}_sT(g\mathbf{y})\rangle&=\langle {}_sT(\mathbf{x}),{}_sT(\mathbf{y})\rangle.
\end{aligned}
\]
The field ${}_sT(\mathbf{x})$ has the form
\begin{equation}\label{sfield}
{}_sT(\theta,\varphi)=\sum_{\ell=s}^{\infty}\sum_{m=-\ell}^{\ell}{}_sa_{\ell m}\,
{}_sY_{\ell m}(\theta,\varphi),
\end{equation}
where ${}_sa_{\ell m}$ are finite variance random variables, that for all  $\ell,\ell'\ge s,$ $ m=-\ell,...,\ell,$  and $m'=-\ell',...,\ell',$
satisfy
\[\mathsf{E}[{}_sa_{\ell m}]=0\quad \mbox{and} \quad
\mathsf{E}[{}_sa_{\ell m}\overline{{}_sa_{\ell'm'}}]=\delta_{\ell\ell'}\delta_{mm'}{}_sC_{\ell},
\]
with ${}_sC_{\ell}\geq 0$ and
\[
\sum_{\ell=s}^{\infty}(2\ell+1){}_sC_{\ell}<\infty.
\]
Note that the series (\ref{sfield}) converges in mean-square in the Hilbert space of square-integrable random sections of the line bundle $(E_s,\pi,\mathbb{S}^2)$, in contrast to the series (\ref{eq:1}) which converges in the Hilbert space of square-integrable random functions on the sphere.

\section{Spin random fields in the ball}\label{sec:ball}

Let us consider a mean-square continuous random field in the ball $\mathbb B(r_0).$ It will be called a spin random field if all its restrictions to  centred spheres of radius $r\in (0, r_0]$ are isotropic spin random fields. In the following the notation ${}_sT(r,\theta,\varphi)$ will be used to denote such fields. Then one obtains
\[
{}_sT(r,\theta,\varphi)=\sum_{\ell=s}^{\infty}\sum_{m=-\ell}^{\ell}{}_sa_{\ell m}(r)\,{}_sY_{\ell m}(\theta,\varphi),
\]
where ${}_sa_{\ell m}(r),$ $r\in[0, r_0],$  are finite variance stochastic processes, that for all  $\ell,\ell'\ge s,$ $ m=-\ell,...,\ell,$   $m'=-\ell',...,\ell',$ and $r,r_1,r_2\in[0, r_0],$   satisfy
\begin{equation}\label{Cl}\mathsf{E}[{}_sa_{\ell m}(r)]=0 \quad \mbox{and} \quad
\mathsf{E}[{}_sa_{\ell m}(r_1)\overline{{}_sa_{\ell'm'}(r_2)}]
=\delta_{\ell\ell'}\delta_{mm'}{}_sC_{\ell}(r_1,r_2),
\end{equation}
with
\begin{equation}\label{Clsum}
\sum_{\ell=s}^{\infty}(2\ell+1){}_sC_{\ell}(r,r)<\infty, \quad r\in [0, r_0].
\end{equation}

To compute the two-point correlation function of the random field ${}_sT(r,\theta,\varphi)$, one can use the addition theorem for spin weighted spherical harmonics. Consider $\mathbf{x}_i\in \mathbb R^3$, $i=1, 2,$ and $\mathbf{e}_z=(0,0,1)^{\top}.$ Let $g_i$ be the rotation with Euler angles $(\varphi_i,\theta_i,0)$ which transforms $e_z$ into $\mathbf{x}_i/\|\mathbf{x}_i\|,$ $i=1, 2$. Let $(\alpha,\beta,\gamma)$ be the Euler angles of the rotation $g^{-1}g_2$. Then
\[
\sum_{m'=-\ell}^{\ell}{}_sY_{\ell m'}(\theta_2,\varphi_2)
\overline{{}_{-m}Y_{\ell m'}(\theta_1,\varphi_1)}
=\sqrt{\frac{2\ell+1}{4\pi}}\,{}_sY_{\ell m}(\beta,\alpha)\mathrm{e}^{-\mathrm{i}s\gamma}.
\]
Using this equation, we obtain
\begin{equation}\label{twopoint}
\langle {}_sT(r_1,\theta_1,\varphi_1),{}_sT(r_2,\theta_2,\varphi_2)\rangle
=\frac{1}{2\sqrt{\pi}}\sum_{\ell=s}^{\infty}\sqrt{2\ell+1}{}_sC_{\ell}(r_1,r_2)
{}_sY_{\ell\,(-s)}(\beta,\alpha)\mathrm{e}^{-\mathrm{i}s\gamma}.
\end{equation}

\begin{remark} Note that the random field ${}_sT(\mathbf{x})$ is mean-square continuous if and only if its two-point correlation function $\langle {}_sT(\mathbf{x}),{}_sT(\mathbf{y})\rangle$ is continuous at all points of the ``diagonal'' set $\{(\mathbf{x},\mathbf{y})\in \mathbb B(r_0)\times\mathbb B(r_0) : \mathbf{x}=\mathbf{y}\}.$ Then, as $|{}_sY_{\ell\,(-s)}(\beta,\alpha)|\le\sqrt{(2\ell+1)/(4\pi)},$ it follows from (\ref{Cl}), (\ref{Clsum}), and (\ref{twopoint}) that to guarantee mean-square continuity each function ${}_sC_{\ell}(r_1,r_2),$  $\ell\ge s,$ must be continuous on  the diagonal set $\{(r_1,r_2)\in [0,r_0]^2: r_1=r_2\}.$
\end{remark}

The stochastic processes ${}_sa_{\ell m}(r),$ $r\in[0, r_0],$ are defined as
\[
{}_sa_{\ell m}(r)=\int_{0}^{\pi}\int_{0}^{2\pi}{}_sT(r,\theta,\varphi)
\overline{{}_sY_{\ell m}(\theta,\varphi)}\sin\theta\,\mathrm{d}\theta
\,\mathrm{d}\varphi.
\]

Let us consider the case when the processes ${}_sa_{\ell m}(r)$ are Gaussian and have continuous sample paths almost surely.  For each $\ell\ge s,$ let  ${}_s\mu_{\ell}$ be the Gaussian probabilistic measure on the Banach space $C([0,r])$ of continuous functions on the interval $[0,r_0]$ that corresponds to the processes ${}_sa_{\ell m}(r).$
By the definition of ${}_sa_{\ell m}(r)$ the measure ${}_s\mu_{\ell}$ is same for all $m=-\ell,...,\ell.$ Let ${}_sH_{\ell}$ be the reproducing kernel Hilbert space of the measure ${}_s\mu_{\ell}$. Finally, let the set $\{\,{}_sf^{(n)}_{\ell}(r)\colon n\in{}_s\mathcal{N}_{\ell}\,\}$ be a Parseval frame in the space ${}_sH_{\ell}$, that is, the set ${}_s\mathcal{N}_{\ell}$ is at most countable, and for any $f\in{}_sH_{\ell}$ it holds
\[
\sum_{n\in{}_s\mathcal{N}_{\ell }}|(f,{}_sf^{(n)}_{\ell})|^2=\|f\|^2,
\]
see \cite{MR3495345}.

 By the result of \cite{MR2511282}, the Gaussian process ${}_sa_{\ell m}(r)$ can be expanded into the series
\begin{equation}\label{eq:4}
{}_sa_{\ell m}(r)=\sum_{n\in{}_s\mathcal{N}_{\ell}}{}_sX^{(n)}_{\ell m}{}_sf^{(n)}_{\ell}(r),
\end{equation}
where ${}_sX^{(n)}_{\ell m}$ are independent standard normal random variables. Moreover, the series \eqref{eq:4} converges uniformly a.s.

In this case
\begin{equation}\label{Clf}
{}_sC_{\ell}(r_1,r_2)=\sum_{n\in{}_s\mathcal{N}_{\ell}}{}_sf^{(n)}_{\ell}(r_1){}_sf^{(n)}_{\ell}(r_2).
\end{equation}

 Conversely, if a stochastic process ${}_sa_{\ell m}(r)$ can be represented in the form of the uniformly a.s. convergent series \eqref{eq:4}, then the set $\{\,{}_sf^{(n)}_{\ell}(r)\colon n\in{}_s\mathcal{N}_{\ell}\,\}$ is a Parseval frame in the space ${}_sH_{\ell}$.

 Finally, by combining the above results, one can see that the random field ${}_sT(r,\theta,\varphi)$ has the following representation
\begin{equation}\label{eq:sprep}
{}_sT(r,\theta,\varphi)=\sum_{\ell=s}^{\infty}\sum_{n\in {}_s\mathcal{N}_{\ell}}\sum_{m=-\ell}^{\ell}{}_sX^{(n)}_{\ell m}\,{}_sf^{(n)}_{\ell }(r)\,{}_sY_{\ell m}(\theta,\varphi).
\end{equation}

See also related wavelet expansions in \cite{10.1093/mnras/stw3176} and \cite{PhysRevD.92.123010}.

\begin{example}

Zernike polynomials in the two-dimen\-sional disk were introduced by \cite{ZERNIKE1934689} to describe aberrations of a lens from the ideal spherical shape.

The 3D Zernike radial polynomials are defined by
\[    R_{n\ell}(r) =\begin{cases}
\sqrt{2n+3} \sum_{k = 0}^{\frac{n - \ell}{2}} ( - 1 )^k \binom{\frac{n - \ell}{2}}{k} \binom{n-k-1+3/2}{\frac{n - \ell}{2}}\,r^{n-2k},\ \mbox{if}\
n -\ell\ \mbox{is even;}\\
0, \ \mbox{if}\ n -\ell\ \mbox{is odd.}
\end{cases}
\]

Note that $R_{n\ell}(r)$ are polynomials of degree~$n$ defined for such $n\ge \ell$ that $n-\ell$ is even. Thus, for a fixed $n\geq 0$, the index $\ell$ takes values $n$, $n-2$, \dots, $n-2\left[\frac{n}{2}\right]$ (where $[\cdot]$ denotes the integer part), i.e. values from $n$ to either $s$ or $s+1.$

In this example we consider the functions
\[
{}_sZ^m_{n\ell}(r,\theta,\varphi)=\tilde{R}_{n\ell}(r)\,{}_sY_{\ell m}(\theta,\varphi),\quad r\in[0, r_0], \ \theta \in[0, \pi],\  \varphi \in [0, 2\pi).
\]

First, let us show how to construct $\{{}_sZ^m_{n\ell},\ n\ge \ell, \,\ell\ge s,\ m=-\ell,...,\ell, \}$ to get a complete orthonormal basis in the space of spin-$s$ functions on  the ball $\mathbb B(r_0).$ Because the spin spherical harmonics are orthonormal on the unit sphere, the polynomials $\{\tilde{R}_{n\ell}(r),\ n\ge \ell\}$  must be orthonormal with the weight function $r^2$ on the interval $[0,r_0].$ The weight function appears due to the Jacobian of the conversion to the spherical coordinates in $\mathbb R^3.$

By the identity (39) in \cite{Mathar2009} any power $r^{\ell+2k},$ $k\in \mathbb N,$ can be represented as a linear combination of $\{R_{(\ell+2i)\ell}(r),i=0,...,k\}.$ Noting that
\[\sum_{k=1}^\infty\frac{\ell+2k+1/2}{(\ell+2k+1/2)^2+1}=+\infty,\]
by the M\"untz theorem, see \cite{Operstein}, one obtains that, for each $\ell,$ the sequence $\{R_{n\ell}(r),\ n\ge \ell\}$ is a basis in $L_2[0,1].$

It is known that, see \cite{Mathar2009},
\begin{equation}\label{orthog}
\int_{0}^{1}r^2R_{n\ell}(r)R_{n'\ell}(r)\,\mathrm{d}r=\delta_{nn'}
\end{equation}
and
\[
R_{n\ell}(r)=\sqrt{2n+3}\,r^{\ell}P^{(0,\ell+1/2)}_{(n-\ell)/2}(2r^2-1),
\]
where $P^{(0,m)}_{k}(\cdot)$ are the Jacobi polynomials \cite[Chapter 22]{AbrSteg}.

By the change of variables $\tilde{r}=r_0r$ in (\ref{orthog}), it follows that in the ball $\mathbb B(r_0)$ it holds

\[
\frac{1}{r_0^3}\int_{0}^{r_0}\tilde{r}^2\,R_{n\ell}\left(\frac{\tilde{r}}{r_0}\right)R_{n'\ell}\left(\frac{\tilde{r}}{r_0}\right)\,\mathrm{d}\tilde{r}=\delta_{nn'}
\]
and one can chose
\[
\tilde{R}_{n\ell}(r)=\frac{\sqrt{2n+3}}{r^{\ell+3/2}_0}r^{\ell}P^{(0,\ell+1/2)}_{(n-\ell)/2}
\left(\frac{2r^2}{r^2_0}-1\right).
\]

Thus, for all $\ell\ge s$ the set
$
\{\,\tilde{R}_{n\ell}(r), n\in {}_s\mathcal{N}_{\ell}\},
$  ${}_s\mathcal{N}_{\ell}=\{n: n \ge \ell,\, n-\ell \ \mbox{is even}\},$ forms a basis in the space of square integrable radial functions on  $\mathbb B(r_0).$ Note that in this case ${}_s\mathcal{N}_{\ell}$ does not depend on $s$ and will be
denoted by $\mathcal{N}_{\ell}.$

If the Hilbert–Schmidt integral operator associated to the kernal ${}_sC_{\ell}(r_1,r_2)$ has the eigenfunctions  $\tilde{R}_{n\ell}(r)$ and eigenvalues $A^{(n)}_{\ell},$ then by Mercer's theorem  the equation (\ref{Clf}) can be rewritten as
\begin{equation}\label{Clf1}
{}_sC_{\ell}(r_1,r_2)=\sum_{n\in\mathcal{N}_{\ell}} A^{(n)}_{\ell}\tilde{R}_{n\ell}(r_1)\tilde{R}_{n\ell}(r_2).
\end{equation}
Then, for each $\ell\ge s,$ the set
$
\{\,\sqrt{A^{(n)}_{\ell}}\tilde{R}_{n\ell}(r), n\in \mathcal{N}_{\ell}\}$   forms a Parseval frame in the space~${}_sH_\ell.$
Thus, the representations (\ref{eq:sprep})  of the corresponding spin random fields in the ball $\mathbb B(r_0)$ has the form
\begin{align}
{}_sT(r,\theta,\varphi)&=\sum_{\ell=s}^{\infty}\sum_{n\in \mathcal{N}_{\ell}}\sum_{m=-\ell}^{\ell}{}_sX^{(n)}_{\ell m}\,\sqrt{A^{(n)}_{\ell}}{}_sZ^m_{n\ell}(r,\theta,\varphi)\notag\\
&=\sum_{n=s}^{\infty}\sum_{k=0}^{\left[\frac{n}{2}\right]} \sum_{m=2k-n}^{n-2k}{}_s{X}^{(n)}_{(n-2k) m}\sqrt{A^{(n)}_{n-2k}}\,\tilde{R}_{n(n-2k)}(r)\,{}_sY_{(n-2k) m}(\theta,\varphi).\notag
\end{align}

\end{example}

\section{Vector $\rho$-stationary random fields in the ball}\label{sec:multi}

This section presents some results on the spectral theory of general $\rho$-stationary vector random fields in the ball. It provides an example of the Mat\'{e}rn random field for a non-Euclidean distance $\rho(\cdot).$   The considered approach is opposite to the one in Sections~\ref{sec:isotr} as a projection of  the ball to a sphere in a higher dimensional space is used.

Let $\rho(\mathbf{x},\mathbf{y})$ denote a distance between points $\mathbf{x},\mathbf{y}\in \mathbb B_0(r_0),$ where $\mathbb{B}_0(r_0)=\{\,\mathbf{x}\in\mathbb{R}^3\colon\|\mathbf{x}\|< r_0\,\}$ is an open ball in $\mathbb R^3.$ Let us consider an isometry $\psi:\mathbb B_0(r_0)\to \mathbb S^3_0(1)$ between the metric spaces $(\mathbb B_0(r_0),\rho)$ and $(\mathbb S_0^3(1),\cos(\gamma)),$ where $\mathbb S_0^3(1)$ is a unit sphere in $\mathbb R^{4}$ with the north pole $(0,0,0, 1)$ removed and $\cos(\gamma)$ is a geodesic distance. Let  $\psi^{(-1)}:\mathbb S_0^3(1)\to \mathbb B_0(r_0)$ denote the inverse mapping for $\psi(\cdot).$

\begin{remark}
As $(\mathbb S_0^3(1),\cos(\gamma))$ is a metric space, then any bijection between $(\mathbb S_0^3(1),\cos(\gamma))$ and $\mathbb B_0(r_0)$ induces a distance in $\mathbb B_0(r_0)$ that can be used as $\rho(\mathbf{x},\mathbf{y}).$ In applications, it is common to consider homeomorphic mappings between these spaces.

Note that there are infinitely many such bijections/homeomorphisms and corresponding distances $\rho(\cdot).$ One of the well-known examples is a composition of the stereographic projection and a mapping of $\mathbb R^3$ onto an open ball.
\end{remark}

Let us consider a vector random field  $\mathbf{T}: \mathbb B_0(r_0)\to \mathbb R^k.$

A zero-mean vector random field  $\mathbf{T}(\mathbf{x})=(T_1(\mathbf{x})),..., T_k(\mathbf{x})),$ $\mathbf{x}\in \mathbb B_0(r_0),$ is called $\rho$-stationary if its covariance matrix $\mathbf{B}(\mathbf{x},\mathbf{y})=\mathsf{E}[\mathbf{T}(\mathbf{x})\otimes\mathbf{T}(\mathbf{y})]$ depends only on the $\rho$-distance between points, i.e.
\[\mathbf{B}(\rho(\mathbf{x},\mathbf{y}))=\mathsf{E}[\mathbf{T}(\mathbf{x})\otimes\mathbf{T}(\mathbf{y})]=\mathsf{E}[\mathbf{T}(\mathbf{x_1})\otimes\mathbf{T}(\mathbf{y_1})]=B(\rho(\mathbf{x_1},\mathbf{y_1})),\] for all $\mathbf{x},\mathbf{x}_1,\mathbf{y},\mathbf{y}_1\in \mathbb B_0(r_0)$ such that $\rho(\mathbf{x},\mathbf{y})=\rho(\mathbf{x}_1,\mathbf{y}_1).$

Let us define a spherical random field $\tilde{\mathbf{T}}(\mathbf{s}),$ $\mathbf{s} \in \mathbb S_0^3(1),$ as $\tilde{\mathbf{T}}(\mathbf{s})= \mathbf{T}(\psi^{(-1)}(\mathbf{s})).$

If $\mathbf{T}(\mathbf{x})$ is $\rho$-stationary, then, due to the isometry of $(\mathbb B_0(r_0),\rho)$ and $(\mathbb S_0^3(1),\cos(\gamma)),$ the random field
$\tilde{\mathbf{T}}(\mathbf{s})$ is isotropic on  $(\mathbb S_0^3(1),\cos(\gamma)).$
Therefore, by \cite[Chapter 1, \S 6]{Yadrenko_book}, the field $\tilde{\mathbf{T}}(\mathbf{s}),$ $\mathbf{s} \in \mathbb S_0^3(1),$ can be represented as
\[\tilde{\mathbf{T}}(\mathbf{s})=\sum_{\ell=0}^{\infty}\sum_{m=1}^{(\ell + 1)^2}
\mathbf{a}_{\ell m}S_{\ell m}(\mathbf{s}),\]
where $S_{\ell m}(\cdot),$ $\ell \in \mathbb N_0,$ $m=1,...,(\ell + 1)^2,$ are spherical harmonics in $\mathbb R^{4}.$

The random coefficients $\mathbf{a}_{\ell m}$ in this spectral representation are defined by
\[\mathbf{a}_{\ell m}=\int_{\mathbb S^3_0(1)}\tilde{\mathbf{T}}(\mathbf{s})
\overline{S_{\ell m}(\mathbf{s})}\,\mathrm{d}\sigma(\mathbf{s}),\]
where $\sigma(\cdot)$ denotes the Lebesgue measure on $\mathbb S^3_0(1).$

Thus, a $\rho$-stationary random field $\mathbf{T}(\mathbf{x})$ can be represented as
\[\mathbf{T}(\mathbf{x})=\tilde{\mathbf{T}}(\psi(\mathbf{x}))=\sum_{\ell=0}^{\infty}\sum_{m=1}^{(\ell + 1)^2}
\mathbf{a}_{\ell m}S_{\ell m}(\psi(\mathbf{x})),\]
\[\mathbf{a}_{\ell m}=\int_{\mathbb S^3_0(1)}{\mathbf{T}}(\psi^{(-1)}(\mathbf{s}))
\overline{S_{\ell m}(\mathbf{s})}\,\mathrm{d}\sigma(\mathbf{s}).\]

If the isometry $\psi(\cdot)$ is also a diffeomorphism with the Jacobian $\mathcal{J}(\cdot),$ then the coefficients $\mathbf{a}_{\ell m}$  can be also computed as
\[\mathbf{a}_{\ell m}=\int_{\mathbb B_0(r_0)}{\mathbf{T}}(\mathbf{x})
\overline{S_{\ell m}(\psi(\mathbf{x}))}\mathcal{J}(\mathbf{x})\,\mathrm{d}\mathbf{x}.\]

These random vector coefficients $\mathbf{a}_{\ell m}$ satisfy the conditions
\[\begin{aligned}
\mathsf{E}[\mathbf{a}_{\ell m}]&=\mathbf{0},\\
\mathsf{E}[\mathbf{a}_{\ell m}\otimes\mathbf{a}_{\ell'm'}]&=\delta_{\ell\ell'}\delta_{mm'}\mathbf{b}_{\ell},
\end{aligned}
\]
with such symmetric nonnegative-definite matrices $\mathbf{b}_{\ell},$ $\ell\in\mathbb N_0,$ that
\[
\sum_{\ell=0}^{\infty}(\ell+1)^2\mathbf{b}_{\ell}<\infty.
\]

Hence, by \cite[Chapter 1, \S 6]{Yadrenko_book} and using the relations between Gegenbauer and  Chebyshev polynomials, see~\cite{AbrSteg}, the two-point correlation function of the vector field $\mathbf{T}(\mathbf{x})$ can be represented as
\[\mathbf{B}(\rho(\mathbf{x},\mathbf{y}))=\langle \mathbf{T}(\mathbf{x}),\mathbf{T}(\mathbf{y})\rangle = \langle \tilde{\mathbf{T}}(\psi(\mathbf{x})),\tilde{\mathbf{T}}(\psi(\mathbf{y}))\rangle =\sum_{\ell=0}^{\infty}\sum_{m=1}^{(\ell + 1)^2}
\mathbf{b}_{\ell}S_{\ell m}(\psi(\mathbf{x}))S_{\ell m}(\psi(\mathbf{y}))\]
and the coefficients $\mathbf{b}_{\ell},$ $\ell\in\mathbb N_0,$ can be computed as
\[\mathbf{b}_{\ell}=\frac{\omega_{3}}{\ell+1}\int_{-1}^{1}\mathbf{B}\left(2\sin\left(\frac{t}{2}\right)\right)\, U_{\ell}(t)\,\sqrt{1-t^2}\,\mathrm{d}t,
\]
where $\omega_{d}={2\pi^{d/2}}/{\Gamma(d/2)}$ and $U_{\ell}(\cdot)$ are the Chebyshev polynomials of the second kind.

By the addition theorem for spherical harmonics the two-point correlation function $\mathbf{B}(\cdot)$ also admits the representation
\[\mathbf{B}(\rho(\mathbf{x},\mathbf{y}))=\frac{1}{\omega_{4}}\sum_{\ell=0}^{\infty}(\ell+1)
{U_{\ell}(\rho(\mathbf{x},\mathbf{y}))}\,\mathbf{b}_{\ell}.\]

\begin{example}
To illustrate  this general approach, let us consider $\psi(\cdot)$ which is a superposition of the stereographic projection and a mapping of $\mathbb R^3$ into an open ball.

The stereographic projection from the north pole $(0,0,0,1)$ acts on spherical points $\mathbf{s}=(s_1,s_2,s_3,s_4)\in\mathbb S_0^3(1)$ as
\[(s_1,s_2,s_3,s_4)\to \left(\frac{s_1}{1-s_4},\frac{s_2}{1-s_4},\frac{s_3}{1-s_4}\right). \]
Its inverse mapping is
\[\mathbf{x}=(x_1,x_2,x_3)\to \left(\frac{2x_1}{1+||\mathbf{x}||^2},\frac{2x_2}{1+||\mathbf{x}||^2},\frac{2x_3}{1+||\mathbf{x}||^2},\frac{||\mathbf{x}||^2-1}{1+||\mathbf{x}||^2}\right).\]
The following homeomorphic mapping from $\mathbf{x}=(x_1,x_2,x_3)\in\mathbb R^3$ to $\mathbb B_0(r_0)$ will be used
\[(x_1,x_2,x_3)\to\left(\frac{2r_0}{\pi}\tan^{-1}\left(x_1\right),\frac{2r_0}{\pi}\tan^{-1}\left(x_2\right),\frac{2r_0}{\pi}\tan^{-1}\left(x_3\right)\right).\]
The superposition of these transformations results in the homeomorphism $\psi(\cdot)$ acting as
\[\psi^{(-1)}(\mathbf{s})=\left(\frac{2r_0}{\pi}\tan^{-1}\left(\frac{s_1}{1-s_4}\right),\frac{2r_0}{\pi}\tan^{-1}\left(\frac{s_2}{1-s_4}\right),\frac{2r_0}{\pi}\tan^{-1}\left(\frac{s_3}{1-s_4}\right)\right)\]
and
\[\psi(\mathbf{x})= \left(\frac{2\tilde{x}_1}{1+||\mathbf{\tilde{x}}||^2},\frac{2\tilde{x}_2}{1+||\mathbf{\tilde{x}}||^2},\frac{2\tilde{x}_3}{1+||\mathbf{\tilde{x}}||^2},\frac{||\mathbf{\tilde{x}}||^2-1}{1+||\mathbf{\tilde{x}}||^2}\right),\]
where $\tilde{x}_i=\tan\left({\pi x_i}/{(2r_0)}\right),$ $i=1,2,3.$

Then, the induced distance $\rho(\cdot)$ on $\mathbb B_0(r_0)$ is
\[\rho(\mathbf{x},\mathbf{y})=C\arccos\left({\frac{4\tilde{\mathbf{x}}^{\top}\tilde{\mathbf{y}}+(1-||\tilde{\mathbf{x}}||^2)(1-||\tilde{\mathbf{y}}||^2)}{(1+||\tilde{\mathbf{x}}||^2)(1+||\tilde{\mathbf{y}}||^2)}}\right),\]
where $C$ is a positive constant and $\tilde{y}_i=\tan\left({\pi y_i}/{(2r_0)}\right),$ $i=1,2,3.$

Let us continue Example~\ref{ex1} and consider the  $\rho$-stationary Mat\'{e}rn random field $T(\mathbf{x}),$ $\mathbf{x}\in \mathbb B_0(r_0),$ with respect to the above distance $\rho(\cdot).$ For simplicity and to be able to visualise numerical results the following computations are presented only for the scalar case, i.e. $k=1.$

The covariance function has the form
\begin{equation}\label{mat2}
\left\langle T(\mathbf{x}),T(\mathbf{y})\right\rangle =\frac{%
2^{1-\nu }\sigma ^{2}}{\Gamma \left( \nu \right) }\left( a\rho(\mathbf{x},\mathbf{y}) \right) ^{{}\nu }K_{{}\nu }\left( a\rho(\mathbf{x},\mathbf{y})
\right) \quad
\end{equation}%
with $\sigma ^{2}>0,$ $a>0$ and $\nu >0.$

\begin{figure}[!htb]
 \captionsetup{width=0.45\textwidth}
	\centering
	\begin{minipage}{0.48\linewidth}
			\includegraphics[width=1\linewidth,height=1\textwidth,trim =5mm 0 15mm 0 ,clip]{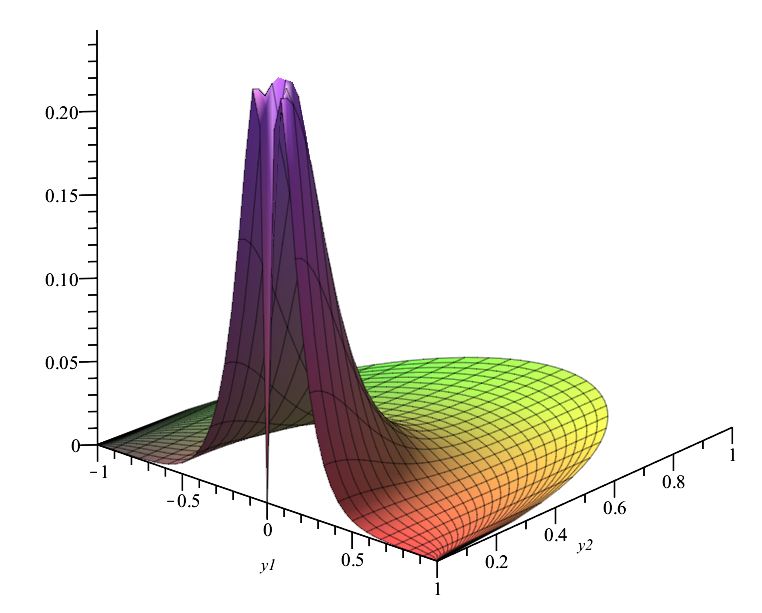}
			\caption{Differences between covariance functions in (\ref{mat1}) and (\ref{mat2}).}\label{fig3}
	\end{minipage}
	\hfill
	\begin{minipage}{0.48\textwidth}
			\includegraphics[width=1\textwidth,height=1\textwidth,trim =5mm 0 5mm 0 ,clip]{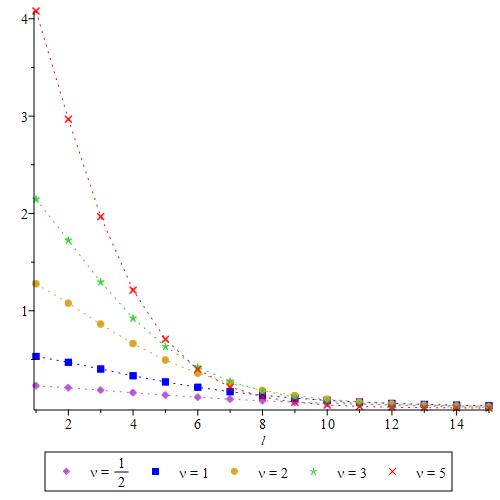}
			\caption{$\mathbf{b}_{\ell}$ for $a=10,$ $\sigma^2=1,$ and various $\nu$.}\label{fig4}
	\end{minipage}
	\end{figure}

The plot of this function is similar to the one in  Figure~\ref{fig1} and is not given here. More informative is Figure~\ref{fig3} which compares this function and the corresponding covariance function from Example~\ref{ex1}, which used the Euclidean distance. To produce the 3D plot the values $\mathbf{x}=\mathbf{0}$ and $\mathbf{y}=(y_1,y_2,y_3)\in \mathbb B_0(r_0)$ with $y_3=0$ were chosen. The horizontal coordinates in Figure~\ref{fig3} are $(y_1,y_2),$ while the vertical one represents the values of the differences between $\left\langle T(\mathbf{0}),T(\mathbf{y})\right\rangle$ in (\ref{mat1}) and (\ref{mat2}).
Figure~\ref{fig3} demonstrates substantial deviations of these two-point correlation functions for distances close to zero.

Because of the isometric mapping, the corresponding covariance function on the sphere $\mathbb S_0^3(1)$ is a restriction of the Mat\'{e}rn stationary covariance function on $\mathbb R^4$ to this unit sphere. Its isotropic spectral density for the 4-dimensional space is
\[f(\lambda)=\sigma
^{2}\,\frac{\nu \left( \nu+1\right) a^{2\nu }}{\left(
a^{2}+\lambda^{2}\right) ^{\nu +2}}.\]

The coefficients $\mathbf{b}_{\ell}$ can be computed by using the formula (22) in \cite[\S 5]{Yadrenko_book} and the result in  Example~\ref{ex1} as
\[
\begin{aligned}
\mathbf{b}_{\ell}&=(2\pi)^4\int_{0}^{\infty}\lambda J_{\ell+1%
}^{2}{(\lambda)} f(\lambda)d\lambda= \sigma
^{2}\,(2\pi)^4\nu \left( \nu+1\right) a^{2\nu }\int_{0}^{\infty}\frac{\lambda J_{\ell+1%
}^{2}{(\lambda)}}{\left(
a^{2}+\lambda^{2}\right) ^{\nu +2}}d\lambda\\
 &=\sigma
^{2}\,\frac{8\pi^4 a^{2\nu }}{\Gamma \left( \nu\right)} \left(
\frac{\Gamma(\ell-\nu)\Gamma(\nu+3/2)}{\sqrt{\pi}\Gamma(\ell+\nu+3)}
{}_1F_2(\nu+3/2;\nu-\ell+1,\nu+\ell+3;a^2)\right.\\
&\quad+\left.\frac{\Gamma(\nu-\ell)a^{2\ell-2\nu}}{2^{2\ell+2}\Gamma(\ell+2)}
{}_1F_2(\ell+3/2;\ell-\nu+1,2\ell+3;a^2)\right).
\end{aligned}
\]

For specific values of the parameters this expressions can be simplified and easily used in computations. For example, for $a=10,$ $\sigma^2=1$ and $\nu=1$ one obtains
\[
\begin{aligned}
\mathbf{b}_{\ell}&=\frac{4\pi^{4}}{25} \left(\Big(\left(l^{2}+3 l +52\right) K_{l +1} \left(10\right)+5 K_{l}\left(10\right) \left(l +2\right)\Big) I_{l +1}\left(10\right)\right.\\
&\quad\left.-5 \Big(\left(l +2\right) K_{l +1}\left(10\right)+10 K_{l}\left(10\right)\Big) I_{l}\left(10\right)\right).
\end{aligned}
\]

For the parameter values $\nu=1/2, 1, 2, 3, 5$ plots of such first spectral coefficients $\mathbf{b}_{\ell}$ are given in Figure~\ref{fig4}.  The plots suggest very fast decay of these coefficients. Thus, in simulations, only few first coefficients can be used to obtain reliable realisations of this $\rho$-stationary Mat\'{e}rn field.

\end{example}

\section{Conclusion}\label{sec:conclusion}
This paper developed the spectral theory for three classes of random fields in the ball. Applications to specific scenarios and the  Mat\'{e}rn correlation model were provided. The derived  spectral representations can be useful for studying theoretical properties and simulating realisations of random fields. Potential areas of    applications include cosmology, geosciences and embryology.

In future studies, it would be also interesting to:
\begin{itemize}
    \item[-] Study rates of convergence in these spectral series representations;
    \item[-] Extend the developed spectral theory to spatio-temporal fields;
    \item[-] Apply the obtained series expansions to investigate evolutions of random fields in the ball driven by SPDEs , see the corresponding results for spherical random fields in \cite{Broadbridge:2019,  Broadbridge:2020, Lang:2015};
    \item[-] Apply the developed methodology to real data, in particular, to new high-resolution cosmological data from future  CMB-S4 and  Euclid mission surveys.
\end{itemize}

\section*{Acknowledgments}

N.~Leonenko and A.~Olenko were partially supported under the Australian Research Council's Discovery Projects funding scheme (project number  DP160101366). We would like to thank Professors Domenico~Marinucci and Ian~Sloan for various discussions about mathematical modelling of CMB data.



\bibliographystyle{amsplain}
\bibliography{Ball}

\end{document}